\documentclass[11pt,twoside]{article}
\usepackage{mathbbold}

\usepackage{amsmath, amssymb, mathrsfs, graphicx}
\usepackage[all]{xy}

\def\scr{\mathscr}

\def\az{\alpha}  \def\bz{\beta}
    \def\dz{\delta}
    \def\fz{\varphi}
  \def\kz{\kappa}
\def\lz{\lambda} 
     \def\oz{\omega}

        \def\sz{\sigma}
        \def\uz{\theta}
\def\vz{\varepsilon}

\def\llz{\Lambda}

\def\qd{\quad}
\def\qqd{\qquad}

\newcommand{\mathsym}[1]{{}}

\def\le{\leqslant}
\def\ge{\geqslant}

\font\cms=cmss9 scaled \magstep1

\def\nnd{\noindent}

\def\thm{\nnd\bg{thm1}}

\def\xmp{\nnd\bg{xmp1}}

\def\dethm{\end{thm1}}

\def\dexmp{\end{xmp1}}

\def\bg{\begin}

\def\de{\end{equation}}

\def\dear{\end{eqnarray}}
\def\lb{\label}
\def\ct{\cite}

\newcommand{\rf}[2]{[\ref{#1}; #2]}

\def\den{\end{enumerate}}

\def\ent{\text{\rm Ent}}
\def\d{\text{\rm d}}
\def\bbb{\mathbbold}

\def\hei{\heiti}

\def\kai{\kaishu}

\begin{document}


\baselineskip 13pt

\thispagestyle{empty}
\renewcommand{\thefootnote}{\fnsymbol{footnote}}

\vspace*{-0.9in}
\noindent {Advances in Mathematics (China) 2017}\newline
\noindent {Vol. 46, No. 4, pp ??}

\begin{center}
{\bf\Large The Charming Leading Eigenpair}
\vskip.15in {Mu-Fa Chen}
\end{center}
\begin{center} (Beijing Normal University)\\
\vskip.1in June 8, 2016       
\end{center}
\vskip.1in

\markboth{\sc Mu-Fa Chen}{\sc The leading eigenpair}



\date{}



\nnd{\bf Abstract}\;\; The leading eigenpair (the couple of eigenvalue and its eigenvector) or the first nontrivial one has different names in different contexts. It is the maximal one in the matrix theory. The talk starts from our new results on computing the maximal eigenpair of matrices. For the unexpected results, our contribution is the efficient initial value for a known algorithm. The initial value comes from our recent theoretic study on the estimation of the leading eigenvalues. To which we have luckily obtained unified estimates which consist of the second part of the talk. In the third part of the talk, the original motivation of the study along this direction is explained in terms of a specific model. The paper is concluded by a brief overview of our study on the leading eigenvalue, or more generally on the speed of various stabilities.
\medskip

\nnd {\small 2000 {\it Mathematics Subject Classification}: 15A18, 65F15, 93E15}

\nnd {\small {\it Key words and phrases}. Leading eigenpair, efficient initial, tridiagonal matrix, speed estimation, Hardy (Poincar\'e)-type inequality, $\fz^4$-model.}

\bigskip

\section{Computing the maximal eigenpair}\lb{s-01}

We begin with the following Perron-Frobenius theorem. For positive $A$ (pointwise), the result is due to Perron, and in the nonnegative irreducible case, it is due to Frobenius. The theorem says there exists uniquely a maximal eigenvalue $\rho(A)>0$ with positive left-eigenvector $u$ and right-eigenvector $g$:
      $$uA = \lz u, \qquad   A g = \lz g, \qquad   \lz= \rho (A).$$
These eigenvectors are also unique up to a constant.

Here is a simplest example due to Luo-Geng Hua (Loo-Keng Hua) (1984)
(refer to \rf{cmf05}{Chapter 10} for references within):

\xmp\lb{t-00}{\rm (Hua, 1984)}\;\;{\cms Let
$$
A =\frac{1}{100}\begin{pmatrix} 25 & 14 \\
40 & 12\end{pmatrix}.$$
Then its maximal eigenvalue $\rho(A)$, the left-eigenvector $u$, and right-eigenvalue
$g$ are, respectively, as follows
$$\aligned
\rho(A)&=\big (37 + \sqrt{2409}\,\big)/200,\\
u&=\big(5\big(13 \!+\! \sqrt{2409}\,\big)/7,\, 20\big)\approx\! (44.34397483, 20),\\
g&=\big(\big(13 \!+\! \sqrt{2409}\,\big)/4,\,20\big).
\endaligned
$$
}\dexmp

Such a simple matrix is already enough to show the great importance of computing the maximal eigenpair.
Recall a simple description of an economic system is using its structure
matrix (the matrix of expanding coefficients) $A$, which is nonnegative, irreducible and invertible. Then the well-known input-output method can be expressed as
$$x_n=x_0 A^{-n},\qquad n\ge 1.$$
where $x_0$ is the input (row vector) and $x_n=\big(x_n^{(0)}, \cdots, x_n^{(d)}\big)$ is the output of the products we are interested at the $n$th year. In 1984, Hua proved the following fundamental theorem:

\thm\lb{t-01}
\emph{\rm (Hua's Fundamental Theorem, 1984)} {\cms
\begin{itemize}    \setlength{\itemsep}{-0.6ex}
\item The optimal choice of $x_0$ is $u$, it has the fastest grow: $x_n=x_0\, \rho(A)^{-n}$.
\item Except some very special $A$, if $x_0\ne u$, then the economic system will be collapsed. That is, some component of the products at some year becomes nonpositive.
\end{itemize}
}\dethm

Certainly, we do not care if the collapse time is very large, say $10^4$ years for instance. However, it is not the case in practice. Table 1 shows the collapse time of Example \ref{t-00} for the initials different from $u$.

\begin{center}{\rm {\bf Table 1} \quad Input and collapse time}\\
\smallskip
\begin{tabular}{|c|c|}
   \hline
${\pmb{x_0^{}}}$ & $\pmb{\text{Collapse time}\;n}$\\
   \hline\hline
 $(44,\; 20)$ & $3$\\
   \hline
$(44.344,\; 20)$ & $8$\\
    \hline
$(44.34397483,\; 20)$ & $13$\\
    \hline
\end{tabular}
\end{center}

If we take only the integer part of $u$ as $x_0$, then the system collapses at the third year; if we take 3 decimals, then the system collapses at the eighth year; finally, if we take all 8  decimals, then the system collapses at the thirteenth year. This result clearly shows the importance of the study on the maximal eigenpair. We need not only high precision but also for large systems.

We now study how to compute the maximal eigenpair.
Before doing so, let us make two remarks.

1) We need to study the right-eigenvector $g$ only. Otherwise, use the transpose $A^*$ instead of $A$.

2) The matrix\,$A$\,is required\,to\,be irreducible with nonnegative off-diagonal elements,
its diagonal elements can be arbitrary. Otherwise, use a shift $A + m I$ for large $m$:
$$(A+mI)g=\lz g\Longleftrightarrow Ag=(\lz-m)g,$$
their eigenvector remains to be the same but the maximal eigenvalues are shifted.

Consider the following example.

\xmp\lb{t-02}{\cms
Consider the matrix
$${Q=\left(
\begin{array}{cccccccc}
 -1 & 1 & 0 & 0 & 0 & 0 & 0 & 0 \\
 1 & -5 & 2^2 & 0 & 0 & 0 & 0 & 0 \\
 0 & 2^2 & -13 & 3^2 & 0 & 0 & 0 & 0 \\
 0 & 0 & 3^2 & -25 & 4^2 & 0 & 0 & 0 \\
 0 & 0 & 0 & 4^2 & -41 & 5^2 & 0 & 0 \\
 0 & 0 & 0 & 0 & 5^2 & -61 & 6^2 & 0 \\
 0 & 0 & 0 & 0 & 0 & 6^2 & -85 & 7^2 \\
 0 & 0 & 0 & 0 & 0 & 0 & 7^2 & -113
\end{array}
\right)}.$$
The main character of the matrix is the sequence $\{k^2\}$.
For this $Q$, the maximal eigenvalue is $-0.525268$ with eigenvector:
$$g\approx (55.878,\; 26.5271,\; 15.7059,\; 9.97983,\;
  6.43129,\; 4.0251,\; 2.2954,\; 1)^*,$$
where the vector $v^*=$ the transpose of $v$.
}\dexmp

Actually, this matrix is truncated from the corresponding infinite one, in which case we have known that the maximal eigenvalue is $-1/4$ (refer to \rf{cmf10}{Example 3.6}).

We now want to practice the standard algorithms in matrix eigenvalue computation. The first method in computing the maximal eigenpair is the
{\it Power Iteration}, introduced in 1929.
Starting from a vector $v_0$ having a nonzero component in the direction of $g$, normalized with respect to a norm $\|\cdot\|$. At the $k$th step, iterate $v_k$ by the formula
$$v_{k}=\frac{A v_{k-1}}{\|A v_{k-1}\|},\quad {z_k}={\|Av_k\|},\qqd k\ge 1.$$
Then we have the convergence: $v_k\to g$ and $z_k\to \rho(Q)$ as $k\to\infty$.
If we rewrite $v_k$ as
$$v_{k}=\frac{A^k v_{0}}{\|A^k v_{0}\|},$$
one sees where the name ``power'' comes from.
For our example, to use the Power Iteration, we adopt the $\ell^1$-norm and
choose $v_0={{\tilde v}_0}/{\|{\tilde v}_0\|}$, where
$${\tilde v}_0\!\!=\!\!(1,\,0.587624,\,0.426178,\,0.329975,\,0.260701,\,0.204394,0.153593,0.101142)^*\!\!.$$
This initial comes from a formula to be given in the last part of this section. Comparing it with $g$,
noting that the eigenvector $g$ decays from 56 to 1, here $\tilde v_0$ decays from 10 to 1, one may worry about the effectiveness of the choice of $v_0$. Anyhow, having the experience of computing its eigensystem, I expect to finish the computation in a few of seconds.
Unfortunately, I got a difficult time to compute the maximal eigenpair for this simple example. Altogether, I computed it for 180 times, not in one day, using 1000 iterations. The printed pdf-file of the outputs has 64 pages. Here are some data.
\def\tba{\begin{array}{cc}
 0 & 2.11289 \\
 1 & 1.42407 \\
 2 & 1.37537 \\
 3 & 1.22712 \\
 4 & 1.1711 \\
 5 & 1.10933 \\
 6 & 1.06711 \\
 7 & 1.02949 \\
 8 & 0.998685 \\
 9 & 0.971749 \\
 10 & 0.948331
  \end{array}}

 \def\tbb{\begin{array}{cc}
  50 & 0.664453 \\
 100 & 0.589332 \\
 200 & 0.542423 \\
 300 & 0.529909 \\
 400 & 0.526517 \\
 500 & 0.525603 \\
 600 & 0.525358 \\
 700 & 0.525292 \\
 800 & 0.525274 \\
 900 & 0.52527 \\
 \ge 990 & 0.525268
  \end{array}}

 \def\ttc{\begin{array}{c}
 \text{Computing}\\
 \text{180 times,} \\
 10^3\,\text{iterations}, \\
 \text{64 pages}. \\
\text{  } \\
 \text{  } \\
(k, -z_k)
   \end{array}}

\vspace{-0.25truecm}

\begin{center}{\rm {\bf Table 2} \quad Outputs $(k, -z_k)$}
\vspace{-0.25truecm}
$$\begin{matrix} \tba\qd& \ttc\;\qd \tbb \end{matrix}$$
\end{center}


{\begin{center}{\includegraphics[width=12.25cm,height=9.0cm]{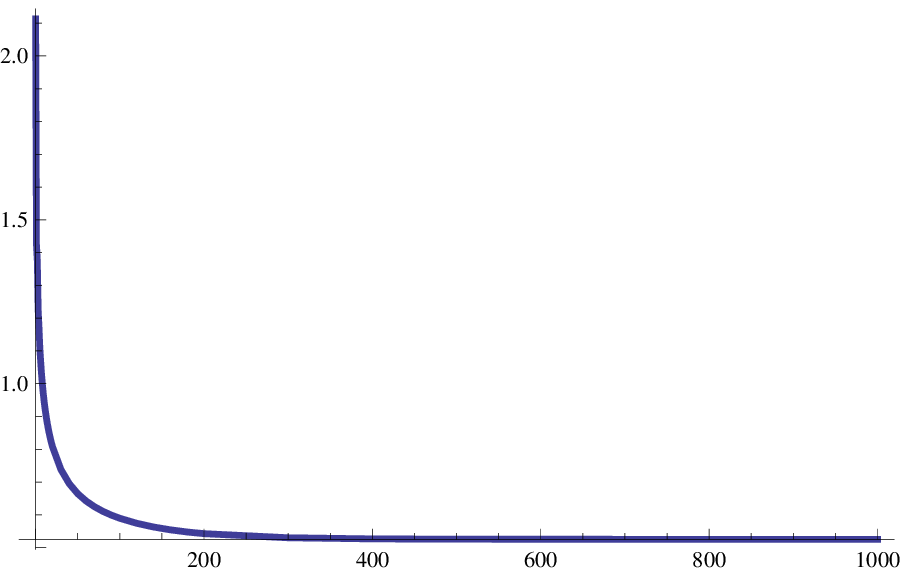}}\newline
{\vspace{-8truecm}\hspace{1truecm}{$$\begin{array}{ll}\text{The figure of } -z_k\\
\text{for }k=0, 1, \ldots, 1000.\end{array}$$}}
\end{center}}
\vspace{6.5truecm}

\nnd The first ten iterations reduce the estimate of the maximal eigenvalue from 2 to 1. It is quite good. Then, we receive the wished output only at the $990$th iteration. The corresponding figure shows that the convergence of $z_k$ goes quickly at the beginning of the iterations.
This means that our initial $v_0$ is good enough. Then the convergence goes very slow which means that the Power Iteration Algorithm converges very slowly.

Let us now move to the second algorithm in computing the maximal eigenpair,
the {\it Rayleigh Quotient Iteration} (RQI), a variant of the {\it Inverse Iteration} presented in 1944.
Here we use the $\ell^2$-norm.
Starting from an approximating pair $(z_0, v_0)$ of the maximal pair $(\rho(A), g)$ with $v_0^*v_0=1$, use the following iteration.
$$ v_k=\frac{(A-z_{k-1}I)^{-1}v_{k-1}}{\|(A-z_{k-1}I)^{-1}v_{k-1}\|}, \qqd z_k=v_{k}^*Av_k,\qqd k\ge 1.$$
If $(z_0, v_0)$ is close enough to $(\rho(A), g)$, then
$$v_k\to g\qd \text{and}\qd z_k\to \rho(A) \qd \text{as }k\to\infty.$$

Before moving further, let us make a remark about this algorithm. Without using the shift $z_{k-1}I$, it is the original
Inverse Iteration:
$$v_k=\frac{A^{-1}v_{k-1}}{\|A^{-1}v_{k-1}\|}\Longleftrightarrow v_k=\frac{A^{-k}v_{0}}{\|A^{-k}v_{0}\|}\qd
\text{i.e. the input-output method}.$$
From this, one may obtain a short proof of Hua's magical assertion in his fundamental theorem. The use of a constant
shift $z I$ for $z$ closed enough to $\rho(A)$ enables us to compute the eigenvector corresponding to $\rho(A)$ rather
than $\lz_{\min}(A)$. The use of a variant shift $z_{k-1}I$ is for accelerating the convergence speed.

Having the hard time spent in the last computation, I was in hesitation to go to the second algorithm. I wondered how many iterations are required using the second algorithm. To have a feeling, I used optimization theory. Suppose we are searching the maximum on the interval (0, 1) for the accuracy of $10^{-6}$. Then, by using the Golden Section Search,
$$   10^{-6} = 0.618^{24}. $$
This means that 24 iterations at least are required. By the Bisection Method,
$$   10^{-6} = 0.5^{20}. $$
Thus, I do not believe that we can complete the job in $20$ iterations.
After prepared enough patient and energy, I started my computation again.
The result came to me, not enough to say surprisingly, I was shocked indeed.

\xmp\;\;{\cms For the same matrix $Q$ and ${\tilde v}_0$ as in Example \ref{t-00}, by RQI, we need two iterations only:}
$$z_1\approx -0.528215,\;\; z_2\approx -0.525268.$$
\dexmp

This shows not only the power of the second method but also the effectiveness of my $v_0$.
For simplicity, from now on, we set $\lz_j:=\lz_j(-Q)$. In particular $\lz_0=-\rho(Q)>0$.

   As usual, ``too good'' is dangerous. For instance, a beautiful person may have a lot of trouble.
Instead of our previous $v_0$, we adopt the uniformly distributed one:
  $$v_0=\{1, 1, 1, 1, 1, 1, 1, 1\}/\sqrt{8}.$$
   This is somehow fair since we may have no knowledge about $g$ in advance.
\xmp\;\;{\cms Let $Q$ be the same as above and use the uniformly distributed $v_0$. Then
$$\aligned
(z_1, z_2, z_3, {\pmb{z_4}}) &\approx  (4.78557,\;5.67061,\;5.91766,\;{\pmb{5.91867}}).\\
(\lz_0, \lz_1, {\pmb{\lz_2}}) &\approx  (0.525268,\,2.00758,\,{\pmb{5.91867}}).
\endaligned$$
}\dexmp
The computation becomes stable at the 4th iteration. Unfortunately, it is not what we want $\lz_0$ but $\lz_2$. In other words, the algorithm converges to a pitfall. Very often, there are $n-1$ pitfalls for a matrix having $n$ eigenvalues. This shows once again our initial ${\tilde v}_0$  is efficient.

   In the last example, $z_0$ is chosen in the automatic way: $z_0=v_0^*(-Q)v_0$.
   If we keep this $v_0$ which is not so good, but using a new $z_0$,
then we come back to our result in two iterations.

\xmp\;\;{\cms Let $Q$ and $v_0$ be the same as in the last example. Choose
$$z_0=2.05768^{-1}\approx 0.485985.$$
Then  $z_1\approx 0.525313,\;\; z_2\approx 0.525268.$
}\dexmp

This shows that the new $z_0$ ($=\dz^{-1}$ to be specified at the end of this section) is efficient.

   We have now computed the same example in 4 times. Here is the comparison of different initials.

\begin{center}{\rm {\bf Table 3} \quad Comparison of different initials}\\
\begin{center}
  \begin{tabular}{| c | c | c |c| }
    \hline
\; {$\pmb Q$}\; &\; ${\pmb{v_0}}$ \;&\; ${\pmb{z_0}}$\; & {\pmb{\# of Iterations}}\\ \hline\hline
   1 & $\pmb{{\tilde v}_0}$ & Power & $10^3$\\ \hline
    2 & $\pmb{{\tilde v}_0}$ & Automatic & {\pmb 2}\\ \hline
    3 &Uniformly distributed & Automatic & Collapse\\ \hline
    4 & Uniformly distributed & $\pmb{\dz_1^{-1}}$ & {\pmb 2}\\ \hline
  \end{tabular}
 \end{center} \end{center}
We now come to the following conclusion.
\begin{itemize}  \setlength{\itemsep}{-0.6ex}
\item RQI is much efficient than Power One.
\item The initials $(v_0, z_0)$ are very sensitive and our $\tilde v_0$ and $z_0=\dz_1^{-1}$ are efficient.
\item It is very hard to handle with the initials. Actually, a large part
  of mathematics research are devoted to this problem.
\end{itemize}

Hopefully, everyone here has heard the name {\it Google's PageRank}. In other words, the Google's search is based on the maximal left-eigenvector
(Exactly the same as what used in the Hua's Theorem \ref{t-01}). On this topic, the following book was published 10 years ago:\\
\text{\qd }Langville, A.N. and Meyer, C. D. (2006).\\
\text{\qd }{\it Google's PageRank and Beyond: The Science of Search Engine Rankings.}\\
\text{\qd }Princeton University Press.\\
In this book, the Power Iteration is included but not the RQI.

Up to now, we have discussed only a small size ($8\times 8 \,(N=7)$) matrix. How about large $N$? In computational mathematics, one often expects the number
of iterations grows in a polynomial way $N^{\az}$ for $\az$ greater or equal to 1. In our efficient  case, since $2=8^{1/3}$, we expect to have $10000^{1/3}=22$ iterations. The next page subverts completely my imagination.

\begin{center}{\rm {\bf Table 4} \quad Comparison of RQI for different $N$}\end{center}
\vspace{-0.6truecm}

\begin{center}
{\begin{tabular}{|c|c|c|c|c|}
   \hline
$\pmb{N+1}$ & $\pmb{z_0}$ & $\pmb{z_1}$ & {$\pmb{z_2=\lz_0}$} & {\bf upper/lower} \\
  \hline\hline
$8$ & 0.523309 & {0.525268} & 0.525268 & $1\!+\!10^{-11}$\\
  \hline
$100$ & 0.387333 & 0.376393 & 0.376383 & $1\!+\!10^{-8}$\\
   \hline
$500$ & 0.349147 & 0.338342 & 0.338329 & $1\!+\!10^{-7}$\\
   \hline
\!$1000$\! & 0.338027 & 0.327254 & 0.32724 & $1\!+\!10^{-7}$\\
   \hline
\!$5000$\! & 0.319895 & 0.30855 & 0.308529 & $1\!+\!10^{-7}$\\
   \hline
\!$7500$\! & 0.316529 & 0.304942 & 0.304918 & $1\!+\!10^{-7}$\\
   \hline
$10^4$ & 0.31437 & 0.302586 & 0.302561 & $1\!+\!10^{-7}$\\
   \hline
\end{tabular}
}\end{center}

Here ${\tilde v}_0$ and $\dz_1$ are computed by our general formulas to be defined very soon below and $$z_0={7}/({8\dz_1})+v_0^*(-Q)v_0/8.$$
We compute the matrices of order $8, 100,\ldots , 10^4$ by using MatLab in a notebook, in no more than 30 seconds, the iterations finish at the second step. This means that the outputs starting from $z_2$ are the same and coincide with $\lz_0$. See the first row for instance, which becomes stable at the first step indeed. We do not believe such a result for some days, so we checked it in different ways. First, since $\lz_0=1/4$ when $N=\infty$, the answers of $\lz_0$ given
in the fourth column are reasonable. More essentially, by using the output $v_2$, we can deduce upper and lower bounds of $\lz_0$ (using \rf{cmf10}{Theorem 2.4\,(3)}), and then the ratio upper/ lower is presented in the last column.
For the first row, by using $v_1$ instead of $v_2$, we also have $1+10^{-7}$. In each case, the algorithm is significant up to 6 digits.

   It is the position to write down the formulas of $\tilde v_0$ and $\dz_1$. Then our initial $z_0$ used in Table 4 is a little modification of $\dz_1^{-1}$: a convex combination of $\dz_1^{-1}$ and $v_0^*(-Q)v_0$.

   Let us consider the tridiagonal matrix. Fix $N\ge 1$ and denote by
       $E=\{0, 1, \ldots, N\}$
the set of indices. By a shift if necessary, we may reduce $A$ to $Q$ with negative diagonals:
$Q=A-m I$, $ m:=\max_{i\in E} \sum_{j\in E} a_{ij},$
$$Q=\!\left(\begin{array}{ccccc}
-(b_0+c_0) & b_0 &0&0 &\cdots \\
a_1 & -(a_1 + b_1+c_1) & b_1 &0 &\cdots \\
0& a_2 & -(a_2 + b_2+c_2) & b_2 &\cdots \\
\vdots &\vdots &\ddots &\ddots &\ddots \\
0& 0  & 0 &\; a_N &-(a_N+c_N))
\end{array}\right)\!.$$

Thus, we have three sequences $\{a_i>0\}$, $\{b_i>0\}$, and $\{c_i\ge 0\}$.
Our main assumption here is that the first two sequences are positive. In order to define our initials, we need three new sequences,
$ \{\mu_k\}$\;(speed measure), $\{h_k\}$, and $\{\fz_k\}$.\footnote{A modification of the algorithm here is presented in \rf{cmf17b}{Apendix \S 4.4}.} The sequence $\{\mu_k\}$ uses $\{a_k\}$ and $\{b_k\}$ only, independent of $\{c_k\}$:
$$\mu_0=1,\;\;\mu_n=\mu_{n-1}\frac{b_{n-1}}{a_n},\qqd 1\le n\le  N.$$
Here and in what follows, our iterations are often of one-step.
Next, we define the sequence $\{h_k\}$:
$$h_0=1,\;\; h_n=h_{n-1}r_{n-1},\qqd 1\le n \le N;$$
here we need another sequence $\{r_k\}$:
$$r_0=1+{c_0}/{b_0},\;\;
r_n=1+\frac{a_n+c_n}{b_n}-\frac{a_n}{b_n r_{n-1}},\qqd 1\le n <N.$$
The boundary of $h$ is defined by
$$h_{N+1}= c_N h_N+ a_N (h_N-h_{N-1}).$$
Note that if $c_k\equiv 0$ for $k<N$, then we do not need the sequence $\{h_k\}$,
simply set $h_k\equiv 1$. Having $\{\mu_k\}$ and $\{h_k\}$ at hand, we can define $\{\fz_k\}$ as follows.
$$\fz_n=\sum_{k= n}^N \frac{1}{h_k h_{k+1}\mu_k b_k}, \qqd 0\le n\le N,\;b_N:=1.$$

We are now ready to define $v_0$ and $\dz_1$ (or $z_0$) using the three new sequences.
$$\aligned
&{{\tilde v}_0(i)\!=\!h_i\sqrt{\fz_i}},\;i\le N; \qqd {v_0}\!=\!{\tilde v}_0/\|{\tilde v}_0\|;\qd \|\cdot\|:=\|\cdot\|_{L^2(\mu)}\\
&{\dz_1}\!=\!\max_{0\le n\le N} \bigg[\sqrt{\fz_n} \sum_{k=0}^n \mu_kh_k^2 \sqrt{\fz_k}
  +\!{\fz_n}^{\!\!\!-1/2}\!\!\sum_{n+1\le j \le N}\!\!\mu_jh_j^2\fz_j^{3/2}\bigg]\!\!=:\!{z_0^{-1}}\!\!.
 \endaligned
$$
Note that $v_0$ and $\dz_1$ are explicitly expressed by these three new sequences. In other words,
we have used three new sequences $ \{\mu_k\}$, $\{h_k\}$, and $\{\fz_k\}$ instead of the original three $\{a_i\}$, $\{b_i\}$, and $\{c_i\}$.

   Finally, the RQI goes as follows.
Solve $w_k$:
$$(-Q-z_{k-1}I) w_k=v_{k-1}, \qqd k\ge 1;$$
and define
$$v_k={w_{k}}/\!{\|w_k}\|,\qqd z_k=(v_{k},\, -Q\, v_k)_{L^2(\mu)}.$$
Then
$$v_k\to g \qd \text{and}\qd z_k\to \lz_0\qqd \text{as } k\to\infty.$$

    Certainly, the next step is going to the general matrix from the tridiagonal one. This is possible once we understand the probabilistic meaning of the sequences $\{\mu_k\}$, $\{h_k\}$, and $\{\fz_k\}$. This work is done in \ct{cmf16b} but omitted here. For more recent progress on this topic,
refer to \ct{cmf17a, cmf17b}.

\section{Unified speed estimation of various stabilities}

We are now going to explain the reason why our initials are efficient. The answer comes from the following result about the unified speed estimation of various stabilities. The result is a short summary of a series of the author's papers published during 2010--2014, starting from 1988.
Refer also to \ct{lzw}.

\thm\lb{t-06}{\rm(Informal\,!)}\;\;{\cms For a tridiagonal matrix $Q$ or a one-dimensional elliptic operator
(order $2$) with or without killing on a finite or infinite interval, in each of {twenty cases}, there exist explicit $\dz$, {$\dz_1$}, $\dz_1'$\,(and then $\dz_n, \dz_n'$,\,recursively)
such that $\dz_n'\uparrow$, $\dz_n\downarrow$ and}
$$
{(4 \dz)^{-1}}\!\!\le \dz_n^{-1}\!\le {\lz_0}\le {\dz_n'}^{-1}\!\!\le {\dz^{-1}}\!, \qqd n\ge 1.$$
{\cms Besides, $1\le {\dz_1'}^{-1}/{\dz_1}^{-1}\le 2$}.
\dethm

The initial $\dz_1$ used in the previous section is taken from here in one specific case. Then the ${\tilde v}_0$
used there was originally used in \rf{cmf10}{\S 3} to deduce $\dz_1$.

Certainly, the notation $\lz_0$ and $\dz_{\#}$ here may be changed case by case. For instance,
for the exponentially ergodic rate (or the exponential decay rate), $\lz_0$ is replaced by
$\az^*$. By \rf{cmf10}{Theorems 1.5 and 7.4} (discrete case) and \rf{cmf12}{Theorem 2.1 and Proposition 6.1}
(continuous case), the rate $\az^*$ coincides with $\lz^{\#}$ to be discussed immediately below and so
the study on $\az^*$ is omitted here.

We now leave the matrix situation and move to differential operators.
First, we consider a special case in parallel to the tridiagonal matrix.
Define the operator
$$L^c= a(x)\frac{\d^2}{\d x^2}+ b(x) \frac{\d}{\d x}-c(x),\qd a(x)> 0, \; c(x)\ge 0$$
on $(0, N)$ with $N\le \infty$. Certainly, by a shift if necessary, one may relax the condition ``$c(x)\ge 0$.''
To study the maximal eigenpair of $L^c$, instead of the
triple $(a, b, c)$ of functions, we introduce three functions $\d \mu/\d x$,
$h$, and $\varphi$ as follows. Let
$$\frac{\d{\mu}}{\d x}\!=\!\frac{e^{C}}{a},\qqd C(x)\!:=\!\int_{0}^x \frac{b}{a},$$
where the Lebesgue measure $\d x $ is omitted in the last integral;
let $h$ be positive $L^c$-harmonic: $L^c h=0$; and let
$$\varphi(x)=\int_0^x \frac{e^{-C}}{h^2}.$$
Having $(\d\mu/\d x, h, \varphi)$ at hand, as in the discrete case, we
can define ${\tilde v}_0$ and $z_0=\dz_1^{-1}$ as follows.
$$\aligned
{\tilde v}_0&=h\sqrt{\varphi},\\
\dz_1&=\sup_{0\le x\le N} \bigg[\sqrt{\fz(x)} \int_{0}^x h^2 \sqrt{\fz}\,\d \mu
+{\fz(x)}^{\!-1/2}\!\!\int_{x}^N\!\!
h^2\fz^{3/2}\,\d\mu\bigg].\endaligned$$

We now go to a more general setup. Consider the space $E=(-M, N)$, $M, N\le \infty$ and the eigenvalue problem:
$$\text{ Eigenequation}: \qd L g=-\lz g, \quad g\ne 0$$
for some differential operator $L$.
Here we use codes `D' and `N' to denote the Dirichlet or Neumann boundary, respectively.

{D}: (Absorbing) Dirichlet boundary,

{N}: (Reflecting) Neumann boundary $g'(-M)= 0$,\\
where $g(-\infty):=\lim_{M\to\infty} g(-M)$. Similarly we have $g'(-\infty)$
and others. Correspondingly, we have four types of eigenvalues.
\begin{itemize}    \setlength{\itemsep}{-0.6ex}
\item $\lz^{\text{\rm NN}}$: Neumann boundaries at $-M$ and $N$.
\item $\lz^{\text{\rm DD}}$: Dirichlet boundaries at $-M$ and $N$.
\item $\lz^{\text{\rm DN}}$: Dirichlet at $-M$ and Neumann at $N$.
\item $\lz^{\text{\rm ND}}$: Neumann at $-M$ and Dirichlet at $N$.
\end{itemize}

Given an elliptic operator $L=L^0$:
$$L=a(x)\frac{\d^2}{\d x^2}+ b(x)\frac{\d}{\d x},$$
define the speed measure $\mu$ and scale measure $\hat\nu$, respectively, as follows
$$\frac{\d{\mu}}{\d x}\!=\!\frac{e^{C}}{a},\qd
\frac{\d{{\hat\nu}}}{\d x}\!=\!e^{-C}\!,\qqd C(x)\!:=\!\int_{\uz}^x \frac{b}{a},$$
where $\uz\in (-M, N)$ is a reference point. Then the leading eigenvalues
$\lz^{\#}$ defined above describe, respectively, the following $L^2(\mu)$-exponential
convergence of the semigroup $\{P_t=e^{t L}\}_{t\ge 0}$:
$$\aligned
&\|P_t f\|\le \|f\|\, e^{-{\lz^{\text{\rm NN}}}\, t},\qd \mu (f)\!:=\!\int_E f \d \mu\!=\!0,\\
&\|P_t f\|\le \|f\|\, e^{-{\lz^{\#}}\, t},
\quad t\ge 0,\; f\in  L^2(\mu),\;\; \text{\rm if \# is not NN}.\endaligned$$
Thus, $\lz^{\text{\rm NN}}$ describes the $L^2$-exponentially ergodic rate
and the other $\lz^{\text{\rm \#}}$ describe the $L^2$-exponential decay rate.

Here is our main result in this part of the talk.

\thm\lb{t-07}{\rm(Chen, 2010)}\;\;{\cms For each \# of 4 cases, we have
the following unified estimates
$${\big(4\kz^{\#}\big)^{-1}\le \lz^{\#}\le  \big(\kz^{\#}\big)^{-1}},$$
where
$$\aligned
\big(\kz^{\text{\rm NN}}\big)^{-1}&= \inf_{x<y}
\big\{{\mu}(-M, x)^{-1} + \mu(y, N)^{-1}\big\}{{\hat\nu}}(x, y)^{-1}\\
\big(\kz^{\text{\rm DD}}\big)^{-1}&=\inf_{x\le y}\big\{{\hat\nu}(-M, x)^{-1}
 + {\hat\nu}(y, N)^{-1}\big\}\mu(x, y)^{-1}\\
\kz^{\text{\rm DN}}&=\sup_{x\in (-M, N)}{\hat\nu}(-M, x)\,\mu(x,N)\\
\kz^{\text{\rm ND}}&= \sup_{x\in (-M, N)} \mu(-M, x)\,{\hat\nu}(x,N)
\endaligned$$
and $\mu(\az, \bz)=\int_{\az}^{\bz} \d \mu$.
In particular, $\lz^{\#}>0$ iff $\kz^{\#}<\infty$.
}\dethm
The beauty of the theorem is displayed in the following aspects.
\begin{itemize}  \setlength{\itemsep}{-0.6ex}
\item Each of the estimates has a universal factor 4.
\item Each constant $\kz^{\text{\#}}$ is expressed by $\mu$ and $\hat\nu$ only.
\item In the expressions of $\kz^{\text{\rm NN}}$ and $\kz^{\text{\rm DD}}$, two boundaries are symmetric.
\item An intrinsic relation between the four constants $\kz^{\text{\#}}$ can be expressed as follows.
$$\renewcommand{\arraystretch}{1.1}
\begin{matrix}
{\begin{array}[c]{ccc}
\kz^{\text{\rm DD}}&\stackrel{\text{Remove }{\hat\nu}(y, N)^{-1}}{\xrightarrow{\hspace*{2.2cm}}} & \kz^{\text{\rm DN}}\\
\Big\updownarrow\scriptstyle{\rm Rule} && \Big\updownarrow\scriptstyle{\rm Rule}\\
\kz^{\text{\rm NN}}&\stackrel{\text{Remove }{\mu}(y, N)^{-1}}{\xrightarrow{\hspace*{2.2cm}}} & \kz^{\text{\rm ND}}
\end{array}}\qd& \qd
{\begin{array}[c]{c}
\text{Rule:}\\
\text{Exchange of codes D and N in $\lz^{\text{\#}}$}\\
\Longleftrightarrow\text{exchange $\mu$ and $\hat\nu$ in $\kz^{\text{\#}}$}
\end{array}
}\end{matrix}$$
\end{itemize}

We remark that the theorem is not as simple as it stands.
In the DN case for instance, it was started by G.H. Hardy in
1920 and completed half a century later  by B. Muckenhoupt et al
around 1970. To obtain the answer in the bilateral cases, one
has to wait for another 40 years until 2010. The proofs in
the last cases use three advanced mathematical tools (the
coupling and distance method, the dual technique, and the capacitary method)
and were completed in five steps (refer to \ct{cmf10}).

There are two ways to generalize the above theorem. The first one
is including the potential term $c$, that is, using $L^c$ instead of $L$.
Again, assume $E=(-M, N),\;M, N\le\infty$. First, we consider the
{\it Poincar\'e-type inequalities}:
$$\lambda_{c}^{\!\#}\,\|f\|_{\mu,\,2}^2\le \|f'\|_{\nu, 2}^2+ {\|c f\|_{\mu, 2}^2},$$
where
$$\nu (\d x)=e^{C(x)}\d x,\qquad {\hat\nu} (\d x)=e^{-C(x)}\d x,$$
and $\|\cdot\|_{\mu, p}=\|\cdot\|_{L^p(\mu)}$.
The inequality becomes equality once $f=g$: the eigenfunction corresponding to $\lambda_{c}^{\!\#}$. This explains the relationship between the inequality and its corresponding eigenvalue.
In particular, when $c\equiv 0$, we return to what we have already studied
above:
$$\sqrt{\lambda^{\#}}\,\|f\|_{\mu,\,2}\le \|f'\|_{\nu, 2}.$$
This leads to the second generalization (generalized to the nonlinear situation): the {\it Hardy-type inequalities}:
$$\|f\|_{\mu,\, q}\le A^{\#} \|f'\|_{\nu,\, p},\qqd p, q\in (1,\infty).$$
We use these inequalities to describe the algebraic convergence $t^{-\az}$ for some $\az>0$. Corresponding to $\nu$ in such a general setup, we have
$${\hat\nu}(\d x)= \exp \bigg[-\frac{C(x)}{p-1}\bigg]\d x$$
which goes back to the previous one when $p=2$.
Finally, we can generalize the left-hand side of the last inequality
to a general normed linear space ${\mathbb B}$:
$$\||f|^q\|_{{\mathbb B}}^{1/q}\le A_{{\mathbb B}}^{\#} \|f'\|_{\nu,\, p}.
$$
A particular use of this class of inequalities is to describe the exponential
convergence in entropy. Note that the entropy functional does not belong
to any $L^q$-space:
$$\|f\|_{L^1(\pi)} \le \ent(f) \le \|f\|_{L^{1+\vz}(\pi)}^{1+\vz}, \qqd \vz>0.$$

The {\it normed linear space} $({\mathbb B}, \|\cdot\|_{{{\mathbb B}}}, \mu)$
here means a subset of Borel measurable functions on $(X, {{\scr X}}, \mu)$ having the following norm
$$\|f\|_{\mathbb B} =\sup_{g\in {\scr G}}\int_X |f|\, g \d \mu,$$
for a given  ${{\scr G}}\subset{{\scr X}}/{{\mathbb R}}_+$.
If we set ${{\scr G}}=L^p\, (p>1)$, then ${{\mathbb B}}=L^{p^*}$: $\frac 1 p+\frac 1 {p^*}=1$.
In the study of logarithmic Sobolev inequality, we use
$${{\scr G}}=\bigg\{g\ge 0: \int_X e^g\d\pi\le e^2+1\bigg\}.$$

Here is a summary of 16 criteria included in Theorem \ref{t-06} (Recall that, as mentioned before, for the omitted 4 cases of $\az^{\text{\#}}$,
we have $\az^{\text{\#}}=\lz^{\text{\#}}$).

\thm\;{\cms The optimal constants $\lz^{\#}$ in the Poincar\'e-type inequalities, with/ without $c$, satisfy
$$\kz^{\#} \le  \lz^{\#\,-1} \le   4\kz^{\#};$$
and the optimal constants $A^{\#}$ in the Hardy-type inequalities, with/without ${\mathbb B}$, satisfy
$$B^{\#} \le  A^{\#} \le  2 B^{\#},$$
where in the DD case for instance, we have

\begin{center}{\rm {\bf Table 5} \quad Isoperimetric constants in different cases}\end{center}\vspace{-0.5truecm}
{\begin{center}{\begin{tabular}{|c|c|}
   \hline
$B_{\mathbb B}$ & $\!\sup\limits_{{x\le y}} \dfrac{\|\bbb{1}_{(x, \, y)} \|_{{\mathbb B}}^{1/q}}
{\big\{{\hat\nu}(-M, x)^{1-p}+ {\hat\nu}(y, N)^{1-p}\big\}^{1/p}}$\\
   \hline
 ${\begin{matrix}
{ {\mathbb B}\!=\!\!L^1(\mu)}\\
 B\end{matrix}}$ & ${\!\sup\limits_{{x\le y}}
\dfrac{{\mu(x, y)^{1/q}}}
{\big\{{{\hat\nu}(-M, x)}^{1-p}+{\hat\nu}(y, N)^{1-p}\big\}^{1/p}}}$\\
    \hline
 ${\begin{matrix}
 {q\!=\!p\!=\!2}\\
 \kz\end{matrix}}$ & ${\!\sup\limits_{{x\le y}}
\dfrac{{\mu(x, y)}}
{{{\hat\nu}(-M, x)}^{-1}+{\hat\nu}(y, N)^{-1}}}$\\
    \hline
 ${\begin{matrix}
 {\text{\rm Killing }c}\\
 \kz_c\end{matrix}}$ & ${\!\sup\limits_{{x\le y}}
\dfrac{{\mu_c}(x, y)}
{{{{\hat\nu}_c}(-M, x)}^{-1}+{{\hat\nu}_c}(y, N)^{-1}}}$\\
    \hline
\end{tabular}
}\end{center}}
\noindent In details, the first line is the most general case
$ {\mathbb B}$. Setting
${\mathbb B}$ to be $L^1(\mu)$, we get the second line, that is the Hardy-type inequalities for $q\ge p$.
Setting $q=p=2$, we get the Poincar\'e-type without $c$. By a change of $\mu$ and $\hat\nu$, we obtain
the last line with $c$: $\mu_c=h^2\mu$, $\hat\nu_c=h^{-2}\hat\nu$, and $h$ is $L^c$-harmonic: $L^c h=0$.
}\dethm

It is remarkable that the previous proofs for the linear case ($q=p=2$) do not suitable to the
present nonlinear situation. To which, we use new analytic proofs
(refer to \ct{cmf13a} and \ct{lzw}).

\section{Original motivation: study on phase transitions}\lb{s-03}

One may be disappointed if I say nothing for the higher dimensional case since up to now we have worked only in dimension one. For this, let us recall the exponential convergence in $L^2$ or in entropy.

Let $\pi$ be a probability measure and denote by $\|\cdot\|$ and $(\cdot, \cdot)$ the norm and inner
product on $L^2(\pi)$. For a given self-adjoint operator $L$ on $L^2(\pi)$:
$$(f, L g)=(Lf, g),\qqd f, g\in {\scr D}(L)\subset L^2(\pi),$$
denote by $\{P_t=e^{t L}\}_{t\ge 0}$ be the semigroup generated by $L$.
We have already seen the {\it exponential stability in $L^2$-sense}:
$$\|P_t f-\pi(f)\| \le  \|f\| e^{-\vz t},\qqd t\ge 0,\; f\in L^2(\pi),$$
and moreover $\vz_{\max}=\lz_1:=\lz^{\text{\rm NN}}.$ Here is an often stronger stability,
{\it exponential stability in entropy}:
$$\aligned
&\ent(P_t f)\le \ent (f)e^{-2{\sz} t},\qqd t\ge 0,\\
&\ent(f)\!:=\!H(\mu\|\pi)\!=\!\!\!\int_E\!\! f \log f \d \pi,\;\;\text{\rm if } \;\frac{\d \mu}{\d \pi}\!=\!f
\endaligned$$

We now go to an infinite-dimensional model. For each $x: {\mathbb Z}^d\to {\mathbb R}$, the interaction potential is
$H(x)=-2J \sum_{\langle i, j\rangle}x_i x_j$ for some $J\ge 0$, where $\langle i, j\rangle$ is the nearest neighbors in ${\mathbb Z}^d$. At each site $i\in {\mathbb Z}^d$, we have the spin potential
$$u(x_i)=x_i^4 - \bz  x_i^2,\qqd x_i\in {\mathbb R},\;\bz\ge 0.$$
The operator for the whole system is
$$L = \sum_{i\in {\mathbb Z}^d} \big[\partial_{ii} - (u'(x_i\!) + \partial_i H)\partial_i \big].$$
\vspace{-0.5truecm}

Here is our main result for this model \big(the $\fz^4$-model\big).

\thm\lb{t-09}{\rm(Chen, 2008)}
$$\begin{aligned}
\inf_{{\Lambda\Subset {\mathbb Z}^d}}\inf_{{\omega\in
{\mathbb R}^{\mathbb Z^d}}} \lambda_1^{\bz,J}\big(\Lambda, \omega \big)
&\approx \inf_{{\Lambda\Subset
{\mathbb Z}^d}}\inf_{{\omega\in {\mathbb R}^{\mathbb Z^d}}} \!\sz^{\bz,J} \!\big(\Lambda, \omega \big)\\
&\approx {\exp\big[-\bz^2/4-c\log \bz\big]}-4dJ\qqd \fbox{$c\!=\!c(\beta)\!\in\! [1, 2]$}\end{aligned}$$
\vspace{-0.8truecm}

{\hspace{-0.5truecm}\includegraphics[width=12.5cm,height=7.5cm]{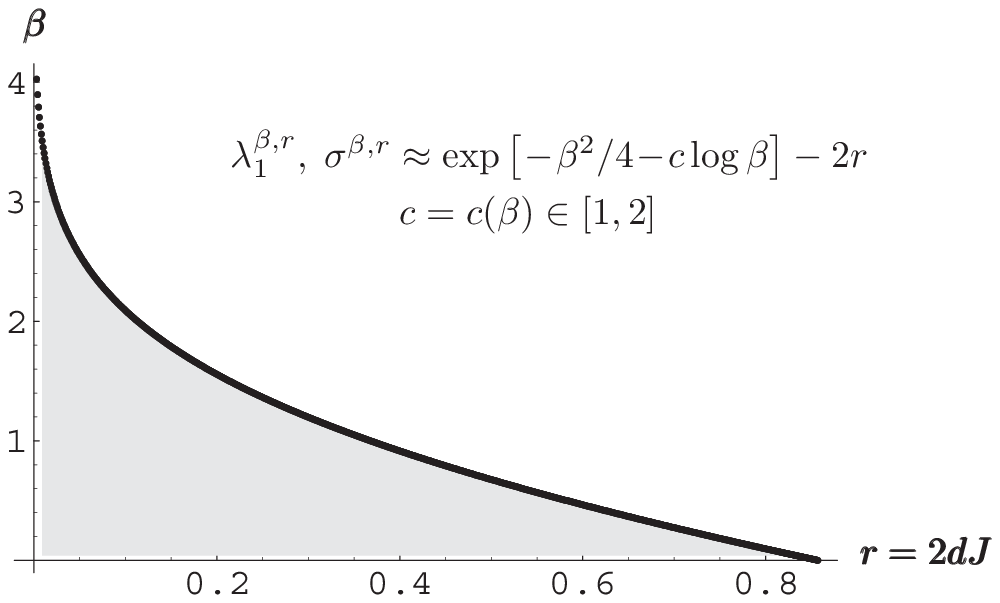}}
\dethm

\nnd Then we proved that the eigenvalue $\lz_1$, as well as the logarithmic Sobolev constant $\sz$ have the same leading decay rate $\exp[-\bz^2/4]-2 r$. More precisely, it says that these two constants have locally such a decay rate uniformly in the finite box $\llz$ and the boundary $\oz$. These constants decay from positive to zero rapidly. This shows the phase transitions of the model. We mention that the leading term $-\bz^2/4$ is exact.

The model illustrates our original motivation of the study on the leading eigenvalue, to describe the phase transitions. Note that for infinite-dimensional mathematics, the known mathematical tools are very limited. We need to look for new mathematical tools. The goal of our study
is developing a new way to describe the phase transitions in statistical
physics. Mathematically, we are looking for a theory of stability speed, an advanced stage of the study on stability. No doubt, such a theory is valuable, as illustrated by Section 1 of the talk.

Up to now, we have discussed the easier part of Theorem \ref{t-06}: $(4\dz)^{-1}\le \lz_0\le \dz^{-1}$, but have not touched the harder
part: $\dz_n^{-1}\le \lz_0\le {\dz_n'}^{-1}$. Hence we have not explained
the way to construct ${\tilde v}_0$ and $\dz_1$ used in \S 1. In the present
situation, we may assume that $h_i\equiv 1$ (otherwise, use \ct{chzh, cmf14} to reduce to this case). Then $\dz_1$ is defined by \rf{cmf10}{(3.4)} and ${\tilde v}_0$ is
the function $f_1$ defined in \rf{cmf10}{Theorem 3.2\,(1)}. Therefore, to
understand $({\tilde v}_0, \dz_1)$, it suffices to have a look at
the first three sections of \ct{cmf10}. We are not going to the details
here. Instead, we prefer to have a short overview of our story, given below.

\section*{Appendix. A brief overview of the research roadmap}

Here we introduce our research roadmap of the topic, and to provide some additional survey articles for the developments of the story.

In 1960's, as a product of the interaction between probability theory and
statistical physics, new branches of mathematics
appeared, first the random fields and then the interacting particle systems,
for instance. We came to the interacting field in 1978, emphasized on the
mathematical foundation of non-equilibrium particle systems.
Our research results were partially collected in \ct{cmf04}. As we know,
a central problem in the study of statistical physics is the
phase transition phenomenon.
Around 1988, we learnt a possible way to describe the phase transition
in terms of the spectral gap (i.e. the first non-trivial eigenvalue,
or more generally the leading eigenvalue) of its generator of the
stochastic process. This led us to a long trip to study the leading
eigenvalue or more generally the speed of various stabilities.

The author's first paper on this topic published in 1991. At the time,
one could compute precisely the principal eigenvalue
of the generator of a Markov Chain in only two or three examples.
This was based on the main theorem in the paper: for a birth--death process,
the ergodic rate (the probabilistic way to describe the
the exponential stability) actually coincides with the first non-trivial eigenvalue
of its generator. If you take a look at this paper and compare it with what I talked above, you will see how far we have come since then. Because our knowledge at the beginning on this topic was rather poor, we started to visit other branches of mathematics.
The first one we visited is the eigenvalue computation for matrices.
In the 1991's paper, we adopted an algorithm to compute the first non-trivial eigenvalue for a class of tridiagonal matrices, without analytic explicit
estimation.

\begin{center}{
{\xymatrix{
\fbox{Statistical Physics}\,\ar@<0.25ex>[r]&\;\fbox{Phase trans 1978}
\ar@<0.5ex>[d]^{1988} \, &\fbox{Probability}\ar@<0.25ex>[l]\\
\fbox{Computational\,Math} \ar@<0.3ex>[r]^{1991}
&\fbox{$\!\!\begin{array}{l}\text{Leading\,eigenvalue}\\\text{Speed of stabilities}\end{array}\!\!$} \ar@<0.3ex>[l]^{{2016}} \ar@<0.5ex>[u]^{{2008}}
& \,\fbox{Riemannian Geom}\;\ar@<0.2ex>[l]^{1998}_{\rm Cheeger}\\
 \ar@<0ex>[ur]_{{2010}}\fbox{3 probability tools}
& \fbox{Hardy ineq} \ar@<1ex>[u]_{{2000}}
&\!\!\fbox{Coupling$+$distance} \ar@<1ex>[u]^{1993}_{{1997}}
  \ar@<0ex>[ul]^{\footnotesize\begin{array}{l}1994\\{1996}\end{array}}_{{1997}}\\
\fbox{Harmonic function} \ar@<1ex>[u]_{{2014}}
& \fbox{Direct proof} \ar@<1ex>[u]_{{2013}}
}}}\end{center}

The next important event is, we found in 1992 that this topic was well studied in
Riemannian geometry. Hence we started to learn the geometric methods, the gradient estimates, in particular. Soon we understood that our probabilistic method --- the coupling method, can also be used for studying this problem. Thus, we went to an opposite way: studying the geometric topic using our probabilistic approach. This was done in several joint papers with Feng-Yu Wang. To obtain sharp estimates however, we need to examine not only the couplings but also the closely related distances.
Thus, the refined method is sometimes called the coupling and distance method. In a survey article of mine, the story was summarized as ``the trilogy of couplings''. The same idea was also used for elliptic operators, as well as matrices. The main credit is that some new variational formula for the lower bound of the eigenvalue was discovered which then improves a number of the known sharp estimates. This may be regarded as our contribution to geometry. After 5 years or so, we also came back to the opposite direction: using some geometric approach (the Cheeger's approach, for instance) to handle with our main problem.

The third important event happened around 2000, we learnt that the Hardy-type inequality (an important subject in Harmonic Analysis) can be used in our study to provide a nice criterion for the positivity
of the principal eigenvalue. This led us to establish 10 criteria for the positive property of different types of stability (or equivalently, inequalities), using our own technique. At the same time, we established new dual variational formulas for the leading eigenvalues, as
well as approximating procedures in computing the eigenvalues.
At this stage, a more or less systematic theory was formed. A series of lectures on the theory up to 2003 consist of the book \ct{cmf05}.

Having worked for 20 years, in 2008, we returned to our original subject, the interacting particle systems (the $\fz^4$-model in particular as discussed in \ref{s-03}) to justify the power of the results obtained until 2003. Luckily, we obtained the exact leading decay rate of the first non-trivial eigenvalue which describes more or less the phase transition curve for the model. We recall that the submission of \ct{cmf08} was delayed for 5 years until we were able to figure out the exact coefficient $1/4$ in the leading rate $\bz^2/4$ given in Theorem \ref{t-09}.

In 2010, we present a unified treatment in \ct{cmf10} of the leading eigenvalue in each of the four cases (i.e. with four different boundary conditions). In this unusually long paper, we obtained not only the unified basic estimates (Theorem \ref{t-07}) but also the improved ones (Theorem \ref{t-06}). Note that the improved estimates are essential for our efficient initials as shown at the beginning of the paper. For this, we have used three probabilistic tools: the coupling and distance method, the dual technique, and the capacitary method. The main ideas of the proofs were surveyed in \ct{cmf12}. Unfortunately, these powerful tools in the linear
case is not suitable for the non-linear one. This is the reason why, to extend the results given in \ct{cmf10} to the Hardy-type inequality, we have to wait for another 13 years. That is, in 2013 (\ct{cmf13a}), we were able to do so by using new direct proofs. Refer to \ct{cmf13b, cmf15} for surveys on \ct{cmf13a}. Thus, only after 13 years known the Hardy-type inequality, we were able to make some contribution to the subject of Hardy inequalities.

The final important event happened in 2014. With Xu Zhang, in \ct{chzhx},
we were able to treat the tridiagonal matrix with general diagonal
elements, using (locally) harmonic functions. This is crucial, otherwise, we can
handle only with a smaller class of tridiagonal matrices (i.e. $c_i\equiv 0$ in the last part of \ref{s-01}). This completes the path $2014\to 2010\to 2016$ in the roadmap above. Recall that we started at using computational mathematics in 1991, and now return to it in 2016, more than 25 years have been passed. All the materials talked here are included in the survey article \ct{cmf16a}
(from which one may find more original references), except part 1 of the talk which has appeared in \ct{cmf16b}.

Sometimes, I feel disappointed since so much time have been spent on
a single topic, I am worrying to be foolish. I tried several times
to leave this area, but I came back, once a new idea appeared, i.e.
the meaning of charming used at the title. Actually, I have been very
lucky for the choice of this topic, so that I can continue my work
for many years, learn much from the other branches of mathematics and
make some contributions to them at last. This overview shows the importance
of choosing a good research topic/direction, and also shows the globality of mathematics. At this moment, I recall that these two points are actually the
main mathematical philosophy presented by D. Hilbert in his famous
lecture given in 1900.

\medskip

\nnd{\bf Acknowledgments}. {\small
This paper is based on a series of talks, five of them given in
2015 were listed in \ct{cmf16b}. The others are presented at Shandon University (2016/3),
Xiamen University (2016/4),
Brigham Young University (2016/4, USA), at the conference
``Frontier Probability Days'' (2016/5, U. of Utah, USA),
at ``The 8th Intern. Conf. on Stoch. Anal. Appl.'' (2016/6, Beijing),
as a distinguished lecture at ``The 10th Cross-Strait Conf. in Stat. \& Probab.'' (abbrev. CSCSP) (2016/8, Chengdu), ``The Sixth National Mathematical Culture Forum'' (2016/8, Lanzhou), Zhejiang University (2016/9), Henan University (2016/9), Southwest Jiaotong University (2016/11),
Center of StatSci, Peking University (2017/3),
and Institute of Mathematics, Academia Sinica (2017/3).
The author acknowledge Professors Shi-Ge Peng, Zeng-Jing Chen,
Ya-Nan Lin, Huo-Xiong Wu, Ke-Ning Lu,
D. Khoshnevisan, E.C. Waymire,
Zhen-Qing Chen, Zhi-Ming Ma,
the organization committee and the local one for CSCSP,
 Jia-An Yan, Hu-Sheng Qiao, Gang Bao, Shu-Xia Feng, Ling-Di Wang, Wei-Ping Li, Shang-Yun Chen,
 Song-Xi Chen, Chii-Ruey Hwang, Tzuu-Shuh Chiang, Yunshyon Chow, Shuenn-Jyi Sheu, and Yuh-Jia Lee
for their invitations and hospitality.
Research supported in part by
         National Natural Science Foundation of China (No. 11131003 and 11626245)
         the ``985'' project from the Ministry of Education in China,
and the Project Funded by the Priority Academic Program Development of
Jiangsu Higher Education Institutions.
}

\vspace{-0.25truecm}

\nnd {\small
Mu-Fa Chen\\
School of Mathematical Sciences, Beijing Normal University,
Laboratory of Mathematics and Complex Systems (Beijing Normal University),
Ministry of Education, Beijing 100875,
    The People's Republic of China.\newline E-mail: mfchen@bnu.edu.cn\newline Home page:
    http://math0.bnu.edu.cn/\~{}chenmf/main$\_$eng.htm
}
\end{document}